\documentclass[13pt]{article}
\usepackage{graphics, graphicx}     

\usepackage{cite}
\usepackage{amsmath}
\usepackage{amssymb}
\usepackage{amsthm}

\theoremstyle{plain}
\newtheorem{Theorem}{Theorem} %
\newtheorem{Lemma}{Lemma}
\newtheorem{Proposition}{Proposition}
\theoremstyle{definition}
\newtheorem{Remark}{Remark}
\newtheorem{Property}{Property}
\newtheorem{Example}{Example}

\theoremstyle{definition}

\newenvironment{Proof} 
{\par\noindent{\it Proof of}} 
{\hfill$\vspace{5mm}\scriptstyle\blacksquare$} 

\numberwithin{figure}{section} 
\numberwithin{table}{section} 

\begin{document}

\setcounter{page}{1}

\markboth{M.I. Isaev, K.V. Isaeva}{On the class of graphs with strong mixing properties}

\title{On the class of graphs with strong mixing properties}
\date{}
\author{ {\bf M.I. Isaev and K.V. Isaeva} }
     
\maketitle
\begin{abstract}
 	We study three mixing properties of a graph:  large algebraic connectivity, large 
 	Cheeger constant (isoperimetric number) and  large spectral gap from 1 for the second largest eigenvalue of the
transition probability matrix of the random walk on the graph. We prove equivalence of this properties (in some sense). We give  estimates for the probability for a random graph to satisfy these properties. In addition, we present  asymptotic formulas for  the numbers of Eulerian orientations and Eulerian circuits in an undirected simple graph.
\end{abstract}

\section{Introduction}
$\ \ \ $
Let $G$ be an undirected simple graph with vertex set $VG$ and edge set $EG$.
	We define $n\times n$ matrix $Q$ by
	\begin{equation}\label{1.1}
		Q_{jk} = 
		\left\{
			\begin{array}{cl}
			-1, & \{v_j,v_k\}\in EG,\\
			\phantom{-}d_j,& j = k,\\
			\phantom{-}0,  & \text{ otherwise,}
		\end{array}\right.
\end{equation}
	where $n = |VG|$ and  $d_j$ denotes the degree of $v_j\in VG$. The matrix $Q = Q(G)$ is called the {\it Laplacian matrix }
	of the graph $G$. The eigenvalues 
	$\lambda_1 \leq \lambda_2 \leq \ldots \leq \lambda_n$ of the
matrix $Q$ are always non-negative real numbers and $\lambda_1 = 0$. The eigenvalue $\lambda_2 = \lambda_2(G)$ is called the
{\it algebraic connectivity} of the graph $G$. The original theory related to algebraic connectivity was produced by Fiedler, see \cite{Fiedler1973}, \cite{Fiedler1989}. The number $\lambda_2(G)$ is a discrete analogue of the smallest positive eigenvalue 
of the Laplace differential operator on the Riemannian manifold. (For more information on the spectral properties of the Laplace matrix see, for example, \cite{Mohar1991} and references therein.)

Let $\cal{F}_\gamma$ be the set  of simple graphs $G$ satisfying  the following property:

\begin{Property}
		\text{The algebraic connectivity} $\lambda_2(G) \geq \gamma |VG|$.
\end{Property}
For a subset of vertices $A \subseteq VG$ let $\partial A$ denote the set of all edges connecting a vertex in $A$ and a vertex outside $A$:
\begin{equation}
	\partial A = \left\{ \{u,v\} \in EG : u\in A, v \in VG \setminus A\right\}.
\end{equation}
The {\it Cheeger constant} (or {\it isoperimetric number}) of $G$, denoted $i(G)$, is defined by
\begin{equation}
	i(G) = \min \left\{ \frac{|\partial A|}{|A|} : A \subset VG,\  0<|A|\leq\frac{|VG|}{2} \right\}.
\end{equation}
The number $i(G)$ is a discrete analogue of the (Cheeger) isoperimetric constant in the theory of Riemannian manifolds and it has many interesting interpretations (for more detailed information see, for example, \cite{Mohar1989} and  references therein).   

Let $\cal{C}_\gamma$ be the set  of simple graphs $G$ satisfying  the following property:

\begin{Property}
		\text{The Cheeger constant (isoperimetric number)} $i(G) \geq \gamma |VG|$.
\end{Property}

Let $P=P(G)$ be the transition probability matrix of the random walk  on  the graph $G$. 
	\begin{equation}\label{1.4}
		P_{jk} = 
		\left\{
			\begin{array}{cl}
			\frac{1}{d_j}, & \{v_j,v_k\}\in EG,\\
			\phantom{-}0,  & \text{ otherwise,}
		\end{array}\right.
\end{equation}
The eigenvalues of $P$ are such that
\begin{equation}
	1 = \chi_1 \geq  \chi_2 \ldots \geq  \chi_n \geq  -1.
\end{equation}
The graph $G$ is connected if and only if the random walk is an irreducible Markov chain. In this case, there exists a unique stationary distribution and multiplicity of the eigenvalue $\chi_1 =1 $ is equal to one. (For more information on random walks on graphs see, for example, \cite{Lovasz1993} and references therein.)

Let $\cal{M}_\gamma$ be the set  of simple graphs $G$ satisfying  the following property:
\begin{Property}
		 The spectral gap $1 - \chi_2(G) \geq \gamma$ and  $\min\limits_j d_j \geq \gamma |VG|$
\end{Property}

Graphs of $\cal{F}_\gamma$, $\cal{C}_\gamma$, $\cal{M}_\gamma$ have strong mixing properties. 
We call $\cal{F}_\gamma \cap \cal{C}_\gamma \cap \cal{M}_\gamma$ as the class of $\gamma$-mixing graphs.
Actually (see Section 2 of the present work), Properties 1-3 are equivalent in the following sense: if a graph satisfies one of these properties with $\gamma = \gamma_0>0$, then it satisfies all  Properties 1-3
with  some $\gamma>0$ depending only on $\gamma_0$. 

In Section 3 we estimate the probability of a random graph to be $\gamma$-mixing. We consider the following model (Gilbert's random graph model):  
every possible edge occurs independently with some fixed probability $0<p<1$. 
It turned out that in this model almost  all graphs (asymptotically) are $\gamma$-mixing with some $\gamma>0$ depending only on $p$.

In Section 4 we construct some general family of graphs, satisfying Properties 1-3  (see Example 3 and Remark 4.1). 
For example, the family of complete bipartite graphs $\{K_{n,n}\}$ is the special case of our general family. We also give some other examples.

In addition, we consider two enumeration problems: counting the number of Eulerian orientations ($EO$) and 
counting the number  of Eulerian circuits ($EC$) in an undirected simple
graph. It is known that both of these problems are complete for the class $\# P$, see 
\cite{Brightwell2005}, \cite{Mihail1996}. 

Recently, in \cite{Isaev2011}, \cite{Isaev2012} the asymptotic behaviour of the
numbers of Eulerian orientations and Eulerian circuits  was determined for $\gamma$-mixing graphs 
(more precisely, for graphs satisfying Property 1). 

In Section 5 we present the asymptotic formulas for $EO$, $EC$ and compare them against the exact
values for small graphs. Actually, if the graph $G$ is $\gamma$-mixing then 
for any $\varepsilon>0$ the error term $|\delta(G)| \leq Cn^{-1/2 + \varepsilon}$, where $C>0$ depends only on $\varepsilon$ and $\gamma$. 
We plan to give the proofs of these formulas in a subsequent paper.     

\section{Equivalence of Properties 1-3}
We recall that for a simple graph $G$ with $n$ vertices: 
\begin{subequations}\label{2.1}
\begin{equation}\label{2.1a}
	\lambda_2(G) \leq \frac{n}{n-1} \min_j d_j,
\end{equation}
\begin{equation}\label{2.1b}
	\lambda_2(G) \geq 2 \min_j d_j - n + 2,
\end{equation}
\end{subequations}
\begin{equation}\label{2.2}
	\frac{\lambda_2(G)}{2} \leq i(G) \leq \sqrt{\lambda_2(G)(2 \max_j d_j - \lambda_2(G))},
\end{equation}

\begin{subequations}\label{2.3}
\begin{equation}
	\lambda_2(G)\leq \lambda_2(G_1)+1,
\end{equation} 
\begin{equation}
	\lambda_2(G)\leq \lambda_2(G'),
\end{equation} 
 where $G_1$ arises from $G$ by removing one vertex  and all adjacent edges, 
 $G'$ is an arbitrary graph such that $VG'=VG$ and $EG \subset EG'$.
\end{subequations}

Estimates \eqref{2.1}, \eqref{2.3}  were obtained in \cite{Fiedler1973}. Estimates \eqref{2.2} were given in \cite{Mohar1989}. 

Using \eqref{2.2} and the inequality $d_j \leq n$, we find that for any $\gamma_0>0$
\begin{equation}\label{FCCF}
	{\cal F}_{\gamma_0} \subset {\cal C}_{\gamma_1} \ \ \text{ and } \ \ {\cal C}_{\gamma_0} \subset {\cal F}_{\gamma_1}, 
\end{equation}
where $\gamma_1>0$ depends only on $\gamma_0$.

In order to complete the proof of the equivalence of Properties 1-3 we need the following lemma. 
For  $\vec{x} \in \mathbb{R}^n$ and $n\times n$ matrix $M$ let us denote
\begin{equation}
	\begin{aligned}
	 \|\vec{x}\| = \sqrt{\vec{x}^T\vec{x}}, \  \ \ \
	 \|M\| = \sup_{\vec{x}\in\mathbb{R}^n,\ \|\vec{x}\|=1} \|M\vec{x}\|.
\end{aligned}
\end{equation}
\begin{Lemma}\label{Lemma1}
	Let $a,b_1,b_2>0$. Let $A$ be  symmetric positive semidefinite $n\times n$ matrix  such that for some 
	$\vec{w} \in \mathbb{R}^n$, $\vec{w}\neq \vec{0}$,
	\begin{subequations}\label{2.6}
		\begin{equation}
			A\vec{w} =0
		\end{equation}
		  and for any $\vec{u}\in \mathbb{R}^n$ such that $\vec{u}^T \vec{w} =0$  
		\begin{equation}
			\|A\vec{u}\| \geq a\| \vec{u} \|.
		\end{equation}
	Then for any  symmetric $n\times n$ matrix $B$ such that 
	\begin{equation}
		\|B\|\leq b_1  \ \ \text{ and } \ \ \vec{w}^T B\vec{w} \geq b_2\|\vec{w}\|^2 
	\end{equation}
	\end{subequations}
	the following statement holds:
	\begin{equation}
		\left\{ 
		\begin{aligned}
			&\det(A - \lambda B) = 0,\\
			&\lambda \neq 0
		\end{aligned}
		\right. 
		\Longrightarrow \lambda \geq \rho
	\end{equation}
	for some $\rho= \rho(a,b_1,b_2) > 0.$
\end{Lemma}

The proof of Lemma \ref{Lemma1} is given at the end of this section.

Note that 
\begin{equation}\label{2.8}
P = I -D^{-1}Q,
\end{equation} 
where $Q$ and $P$ are the same as in \eqref{1.1} and \eqref{1.4}, respectively, $I$ denotes the identity matrix and $D$ is the diagonal 
 matrix defined by $D_{jj} = d_j$.

Let $A_1= \frac{1}{n} Q$, $B_1= \frac{1}{n}D$ and $\vec{w}_1 =[1,\ldots,1]^T$. Using \eqref{2.8}, we find that:

\begin{equation}\label{2.9}
	\det (A_1 - \lambda I) = 0 \Longleftrightarrow \det (Q - \lambda n I) = 0;
\end{equation}
\begin{equation}
	\det (A_1 - \lambda B_1) = 0 \Longleftrightarrow \det(P - (1-\lambda)I)=0;
\end{equation}
for any $\vec{u}\in \mathbb{R}^n$ such that $\vec{u}^T \vec{w}_1 =0$  
\begin{equation}
\vec{u}^T A_1 \vec{u} \geq \frac{1}{n}\lambda_2(G) \vec{u}^T\vec{u};
\end{equation}
\begin{equation}\label{2.12}
	\|B_1\|\leq \frac{1}{n} \max\limits_{j}d_j \leq 1, \ \ \ \vec{w}_1^{T} B_1 \vec{w}_1 \geq \frac{1}{n} \min\limits_{j}d_j \,\vec{w}_1^{T}\vec{w}_1. 
\end{equation}
Combining Property 1, \eqref{2.1a}, \eqref{2.9}-\eqref{2.12} and Lemma \ref{Lemma1}, we get that for any $\gamma_0>0$
\begin{equation}\label{FM}
	{\cal F}_{\gamma_0} \subset {\cal M}_{\gamma_2},
\end{equation}
where $\gamma_2>0$ depends only on $\gamma_0$.

Let $A_2= D^{-\frac{1}{2}} Q D^{-\frac{1}{2}}$, $B= nD^{-1}$ and $\vec{w}_2 =D^{\frac{1}{2}}\vec{w}_1$, where 
$D^s$ is the diagonal 
 matrix defined by $D_{jj}^{s} = (d_j)^s$. Using \eqref{2.8}, we find that:

\begin{equation}\label{2.14}
	\det (A_2 - \lambda I) = 0 \Longleftrightarrow  \det(P - (1-\lambda)I)=0;
\end{equation}
\begin{equation}
	\det (A_2 - \lambda B_2) = 0 \Longleftrightarrow \det (Q - \lambda n I) = 0;
\end{equation}
for any $\vec{u}\in \mathbb{R}^n$ such that $\vec{u}^T \vec{w}_2 =0$  
\begin{equation}
\vec{u}^T A_2 \vec{u} \geq (1 - \chi_2(G)) \vec{u}^T\vec{u};
\end{equation}
\begin{equation}\label{2.17}
	\|B_2\|\leq \frac{n}{ \min\limits_{j}d_j}, \ \ \ \vec{w}_2^{T} B_1 \vec{w}_2 = n \vec{w}_1^{T}\vec{w}_1
	 \geq \vec{w}_2^{T}\vec{w}_2  . 
\end{equation}
Combining Property 3, \eqref{2.14}-\eqref{2.17} and Lemma \ref{Lemma1}, we get that for any $\gamma_0>0$
\begin{equation}\label{MF}
	{\cal M}_{\gamma_0} \subset {\cal F}_{\gamma_3},
\end{equation}
where $\gamma_3>0$ depends only on $\gamma_0$.

Putting together \eqref{FCCF}, \eqref{FM} and \eqref{MF}, we obtain the desired assertion:
\begin{Theorem}\label{Theorem1}
	Let ${\cal F}_{\gamma}, {\cal C}_{\gamma}, {\cal M}_{\gamma}$ be defined as in Section 1. Then for any $\gamma_0>0$
	\begin{equation}
		{\cal F}_{\gamma_0}\cup{\cal C}_{\gamma_0}\cup{\cal M}_{\gamma_0} \subset 
		{\cal F}_{\gamma}\cap{\cal C}_{\gamma}\cap{\cal M}_{\gamma},
	\end{equation}
where $\gamma>0$ depends only on $\gamma_0$.	
\end{Theorem}

Now it remains to prove Lemma \ref{Lemma1}.
\begin{Proof} {\it Lemma \ref{Lemma1}.}
	Let $\det(A-\lambda B) =0$. Then for some $\vec{v} \in \mathbb{R}^n$, $\vec{v}\neq 0$, 
	\begin{equation}\label{2.20}
		A \vec{v} =  \lambda B \vec{v}.
	\end{equation}
Let $\vec{v} = \vec{v}_{\parallel} + \vec{v}_{\bot}$, where $\vec{v}_{\parallel} \parallel \vec{w}$ and $\vec{v}_{\bot}^T \vec{w} =0$.
Due to (\ref{2.6}a), we have that
\begin{equation}
 \vec{v}_{\parallel}^T A \vec{v} = 0. 
\end{equation}
Since $\lambda \neq 0$, using \eqref{2.20}, we get that
\begin{equation}\label{2.22}
	\vec{v}_{\parallel}^T B \vec{v}_{\parallel} = -  \vec{v}_{\parallel}^T B \vec{v}_\bot.
\end{equation}
Using (\ref{2.6}c), \eqref{2.22} and the Cauchy–Schwarz inequality, we find that
\begin{align}
	b_1\|\vec{v}_{\parallel}\|\|\vec{v}_{\bot}\| \geq \|\vec{v}_{\parallel}\| \|B\vec{v}_{\bot}\| 
	\geq |\vec{v}_{\parallel}^T B \vec{v}_\bot| = |\vec{v}_{\parallel}^T B \vec{v}_{\parallel}| \geq 
	b_2\|\vec{v}_{\parallel}\|^2.
\end{align} 
Thus we have that
\begin{equation}\label{2.24}
	\|\vec{v}_{\bot}\| \geq \frac{b_2}{\sqrt{b_1^2 + b_2^2}} \|\vec{v}\|.
\end{equation}
Using (\ref{2.6}b), (\ref{2.6}c) and \eqref{2.24}, we find that:
\begin{equation}\label{2.25}
	\begin{aligned}
	&\|A\vec{v}\| \geq a \|\vec{v}_{\bot}\|,\\
	&\|B\vec{v}\| \leq b_1 	\|\vec{v}\| \leq \frac{b_1\sqrt{b_1^2 + b_2^2}}{b_2} \|\vec{v}_{\bot}\|.
	\end{aligned}
\end{equation}
Combining (\ref{2.20}) and \eqref{2.25}, we obtain that
\begin{equation}
	\lambda \geq \frac{ab_2}{b_1\sqrt{b_1^2 + b_2^2}}. 
\end{equation}
\end{Proof}
\section{Probability for a random graph to be $\gamma$-mixing}
Let $\xi$ be a random variable belonging to the binomial distribution $B(M,p)$: 
\begin{equation}
	\mbox{Pr}(\xi = k) = \frac{M!}{k! (M-k)!} p^k (1-p)^{M-k}, \ \ \  0<p<1, \ M\in \mathbb{N}.
\end{equation}
Note that 
\begin{equation}\label{3.2}
	\mbox{Pr}(\xi \leq \alpha M) \leq c^{-M}
\end{equation}
for some $\alpha>0$, $c>1$ depending only on $p$.  This follows, for example, from  the following estimate: 
for $1 \leq k \leq \frac{p(M+1)}{p+2}$
\begin{equation}
	\frac{\mbox{Pr}(\xi = k)}{\mbox{Pr}(\xi = k-1)} = \frac{M-k+1}{k} \, \frac{p}{1-p} \geq 
	\frac{\frac{2}{p+2}(M+1)}{\frac{p}{p+2}(M+1)}\, \frac{p}{1-p} \geq  2.
\end{equation}

Let $G$ be a random graph belonging to Gilbert's random graph model $G(n,p)$:
\begin{equation}
\begin{aligned}
	\forall_{{1 \leq i<j\leq n}} \ \mbox{Pr} (\{v_i,v_j\} \in EG) =  p, \ \ \ 0<p<1,\\  
	\text{(independently for each $\{i,j\}$)}.
\end{aligned}
\end{equation}
For a subset of vertices $A \subset VG$, using \eqref{3.2}, we find that  
\begin{equation}\label{3.5}
	\mbox{Pr}(|\partial A| \leq \alpha |A|(n-|A|)) \leq c^{-|A|(n-|A|)}.
\end{equation}
Using \eqref{3.5}, we get that
\begin{equation}\label{3.6}
	\begin{aligned}
		\mbox{Pr}(i(G) \leq \frac{\alpha n}{2}) 
		&\leq \sum\limits_{A \subset VG,\ 0<|A|\leq\frac{n}{2} } 
		\mbox{Pr}(|\partial A|  \leq \alpha  |A| \frac{n}{2}) \leq \\ 
		&\leq
		\sum\limits_{k=1}\limits^{n/2}\sum\limits_{|A| = k} 
		\mbox{Pr}(|\partial A|  \leq \alpha  |A| (n-|A|))\leq\\ 
		&\leq \sum\limits_{k=1}\limits^{n/2} \frac{n!}{k! (n-k)!} \, c^{-k(n-k)}
		\leq \sum\limits_{k=1}\limits^{n/2} \frac{n!}{k! (n-k)!}\, (c^{-n/2})^k\leq \\
		&\leq (1 + c^{-n/2})^n - 1 \leq \beta^{-n}
	\end{aligned}
\end{equation}
for some $\beta > 1$ depending only on $p$. 

Due to \eqref{3.6} and Theorem \ref{Theorem1}, we obtain that probability for a random graph (in 
Gilbert's model $G(n,p)$) to be 
$\gamma$-mixing is  at least $1-\beta^{-n}$, where $\gamma=\gamma(p)>0$ and $\beta=\beta(p)>1$. 
\section{Some basic properties and examples}
We note that, due to \eqref{2.1b} and Theorem \ref{Theorem1}, 
\begin{equation}
	\begin{aligned}
	\text{ if } &\min\limits_{j}\, d_j \geq \sigma |VG| \\ &\text{ for some } \sigma > 1/2 
	\end{aligned}
	\Longrightarrow 
	\begin{aligned}
	\text{the graph $G$ is } \gamma&\text{-mixing}\\
	\text{ for some } &\gamma= \gamma(\sigma) > 0.
	\end{aligned}
\end{equation}
\begin{Example}
Let $K_n$ and $\tilde{K}_n$ be two complete graphs with $n$ vertices. We define $G^{(1)}_n$ by
\begin{equation}
	\begin{aligned}
	VG^{(1)}_n = VK_n &\cup V\tilde{K}_n,\\ EG^{(1)}_n = EK_n &\cup E\tilde{K}_n \cup E^+, \\
	&\text{ where } E^+  = \{\{v_i,\tilde{v}_i\}, i=1,\ldots n\}.
	\end{aligned}
\end{equation}
\end{Example}
\noindent
For $G=G^{(1)}_n$ we have that for all $j=1, \ldots ,2n$
\begin{equation} 
 d_j = n+1 > \frac{1}{2}|VG^{(1)}_n|,
\end{equation} 
but
\begin{equation}
	 \frac{i(G^{(1)}_n)}{n} \leq \frac{|E^+|}{n|VK_n|} \rightarrow 0, \ \ \text{ as $n\rightarrow \infty$}.
\end{equation}
Thus the family $\{G^{(1)}_n\}$ does not satisfy Property 2 (and hence Properties 1,3). 
We note also that even the vertex and the edge connectivity is large for this family of graphs.
 
\begin{Example}
Let $K_n$  be the complete graph with $n$ vertices. We define $G^{(2)}_n$ by
\begin{equation}
	\begin{aligned}
	VG^{(2)}_n = VK_n &\cup v_{n+1},\\ EG^{(2)}_n = EK_n &\cup \{v_n, v_{n+1}\}.
	\end{aligned}
\end{equation}
\end{Example}

\noindent
For $G=G^{(2)}_n$ we have that $\min\limits_j d_j = 1$, but one can show that $\chi_2(G^{(2)}_n) = 1/\sqrt{n}$. Note also that Examples 1 and 2 show, in particular, that both conditions of Property 3 are necessary.

\begin{Example}
	Let $G^0$ be a connected simple graph with $m>1$ vertices. Let $c_1,c_2,\ldots, c_m$ be some natural numbers.
	We define $G_n^{(3)}$ by
	\begin{equation}
		\begin{aligned}
			VG^{(3)}_n = \{v_j^i:\  i =1, \ldots, nc_j, \ j =1,\ldots, m\},\\ 
			\{v^{i_1}_{j_1}, v^{i_2}_{j_2}\} \in EG^{(3)}_n \Longleftrightarrow 	\{v_{j_1}, v_{j_2}\} \in EG^0.
		\end{aligned}
	\end{equation}
\end{Example}

\noindent
We estimate the Cheeger constant (isoperimetric number) $i(G^{(3)}_n)$. Let 
\begin{equation}
c_{0} = \min\limits_{1\leq j\leq m} c_j, \ \ C = \sum\limits_{j=1}\limits^{m} c_j \ . 
\end{equation} 
We have that 
\begin{subequations}\label{4.8}
\begin{equation}
 \text{ the degree of  each vertex of $VG^{(3)}_n$ at least $c_{0} n$ } 
\end{equation}
 and 
 \begin{equation}
 		|VG^{(3)}_n| = Cn.
 \end{equation}
\end{subequations}

Let $A \subset VG^{(3)}_n$, 
$|A| \leq Cn/2$. 
 \begin{itemize}
 	\item Case 1. $|A|\leq c_0 n/2$. Using \eqref{4.8}, we find that
 	\begin{equation}\label{4.9}
 		|\delta A| \geq c_0 n |A| - |A|^2 \geq  \frac{c_0n}{2} |A| = \frac{c_0|VG^{(3)}_n|}{2C}|A|
 	\end{equation}
 	\item Case 2. $|A|> c_0 n/2$. Let $V_1,V_2 \subset VG_0$ such that
 	\begin{equation}
 		v_j \in V_1 \Longleftrightarrow |\left\{v_j^{i}: v_j^{i} \in A  \right\}| \geq  \frac{c_0n}{2m}, 
 	\end{equation}
 	\begin{equation}
 		v_j \in V_2 \Longleftrightarrow |\left\{v_j^{i}: v_j^{i} \notin A  \right\}| \geq  \frac{c_jn}{2}. 
 	\end{equation}
 	Due to $c_0 n/2 < |A| \leq Cn/2$, we have that
 	\begin{equation}
 			V_1 \cup V_2 = VG_0, \ \ |V_1|>0, \ \ |V_2|>0.
 	\end{equation}
 	Since  $G^0$ is connected, we can find $v_{j_1} \in V_1$, $v_{j_2} \in V_2$ such that $\{v_{j_1}, v_{j_2}\} \in EG_0$. 
 	Estimating the number of edges in $\partial A$, matching these vertices, we obtain that
 	\begin{equation}\label{4.13}
 		|\delta A| \geq \frac{c_0c_{j_2}n^2}{4m} \geq \frac{c_0^2n}{2mC} |A| = \frac{c_0^2	|VG^{(3)}_n|}{2mC^2} |A|.
 	\end{equation}
 \end{itemize}

Combining \eqref{4.9}, \eqref{4.13} and  Theorem \ref{Theorem1}, we obtain that 
	the family $\left\{ VG^{(3)}_n \right\}$ satisfies Properties 1-3 with $\gamma>0$, 
	depending only on $G_0$, $c_1,\ldots,c_m$.
	
We note that Example 3 can be modified so that the constants $c_1,\ldots,c_m>0$   are not necessarily natural numbers. We assumed that just for simplicity of the proof.

\begin{Remark}
Combining  \eqref{2.3}, Theorem \ref{Theorem1} and Example 3, one can prove Properties 1-3
for a large number of classic examples (including $\{K_n\}$, $\{K_{n,n}\}$ and many others).
\end{Remark}
\section{Asymptotic estimates for $\gamma$-mixing graphs}
	An {\it Eulerian orientation} of $G$  is an orientation of its edges such that for every vertex the number of incoming edges and outgoing edges are equal. 
	We denote $EO(G)$ the number of Eulerian orientations. 
	 Eulerian orientations of the complete graph $K_n$ are called regular tournaments. 

 An {\it Eulerian
circuit } in $G$ is a closed walk 
which uses every edge of $G$ exactly once. Let
$EC(G)$ denote the number of these up to cyclic equivalence. 

We recall that $EO(G) = EC(G) = 0$, if the degree of at least one vertex of $G$ is odd (for more information see, for example, \cite{Biggs1976}).	In this section we always assume that every vertex has even degree.

Let consider two enumeration problems: counting the number of Eulerian orientations  and 
counting the number  of Eulerian circuits in an undirected simple
graph. It is known that both of these problems is complete for the class $\# P$, see 
\cite{Brightwell2005}, \cite{Mihail1996}. 

The results presented in this section are based on estimates of \cite{Isaev2011}, \cite{Isaev2012}. 
We plan to give detailed proofs in a subsequent paper.

\subsection{Eulerian orientations}

We recall that the problem of counting Eulerian orientations can be reduced to counting perfect matching for some class of bipartitute graphs for which it can be done approximately with
high probability in polynomial time, see \cite{Mihail1996}. However, the degree of the polynomial is large, so, in fact, these algorithms have a very big work time  for the error term $O(n^{-1/2})$.

For $\gamma$-mixing graphs we have the following asymptotic formula:

\begin{Proposition}\label{main}
	Let $G$ be an undirected simple graph  with  $n$ vertices $v_1, v_2, \ldots, v_n$ having even degrees. Let $G$ be a $\gamma$-mixing graph 
	for some $\gamma>0$. 
	Then 
	\begin{equation}\label{main_eq}
		\begin{aligned}
		EO(G) = \left(1 + \delta(G) \right)
		\left(
		  2^{|EG|+\frac{n-1}{2}} \pi^{-\frac{n-1}{2}} \frac{1}{ \sqrt{t(G)}} \prod\limits_{\{v_j,v_k\}\in EG} P_{jk} 
		\right),\\
			P_{jk} = 1 - \frac{1}{4(d_j+1)^2} -  \frac{1}{2(d_j+1)(d_k+1)} - \frac{1}{4(d_k+1)^2},
			\end{aligned}
	\end{equation}
	where  $d_j$ denotes the degree of vertex $v_j$,	  $t(G)$ denotes the number of spanning trees
of the graph $G$ and for any  $\varepsilon>0$ 
	\begin{equation}
		|\delta(G)| \leq C n^{-1/2+\varepsilon},
	\end{equation}
	where constant $C > 0$ depends only on $\gamma$ and $\varepsilon$.   
\end{Proposition}

\begin{Remark}
We note that, according to the Kirchhoff Matrix-Tree-Theorem, see \cite{Kirchoff1847}, we have that
\begin{equation}\label{Eq_1_3}
	t(G) = \frac{1}{n}\lambda_2\lambda_3\cdots\lambda_{n} = \det{M_{11}},
\end{equation}
where $M_{11}$ results from deleting the first row and the first column of 
	$Q$.
\end{Remark}

\begin{Remark}
	For the complete graph $\lambda_2(K_n) = n$, $EK_n = \frac{n(n-1)}{2}$, $t(K_n) = n^{n-2}$, 
\begin{equation}
	\begin{aligned}
	\prod\limits_{\{v_j,v_k\}\in EK_n} P_{jk} =
	\left(1 - \frac{1}{4n^2} - \frac{1}{2n^2} - \frac{1}{4n^2} \right)^{\frac{n(n-1)}{2}} = \\ =
	\left( e^{\ln\left(1-\frac{1}{n^2}\right)} \right) ^{\frac{n(n-1)}{2}} = e^{-1/2} + O(n^{-1}).
	\end{aligned}
\end{equation}
	The result of Proposition \ref{main} for this case is reduced to the result of \cite{Brendan1990} 
	on counting regular tournaments in the complete graph.
\end{Remark}
We plan to give the proof of Proposition \ref{main} in a subsequent paper. 
In the present work we just compare the answers given by  formula \eqref{main_eq} against the exact
values for small graphs. 
Let
\begin{equation}
	Error(G) = \frac{EO_{approx}(G) - EO(G)}{EO(G)},
\end{equation} 
where $EO_{approx}(G)$ is taken according to right-hand side of \eqref{main_eq}.
The following charts show the dependence of $Error(G)$ on the ratio $\lambda_2(G)/n$, where $\lambda_2(G)$ 
is the algebraic connectivity and $n=|VG|=6,7,8,9$:

\vspace{5mm}

\renewcommand{\baselinestretch}{0.9}
\begin{figure}[h]
\vspace{0mm}
\includegraphics[width=0.45\textwidth]{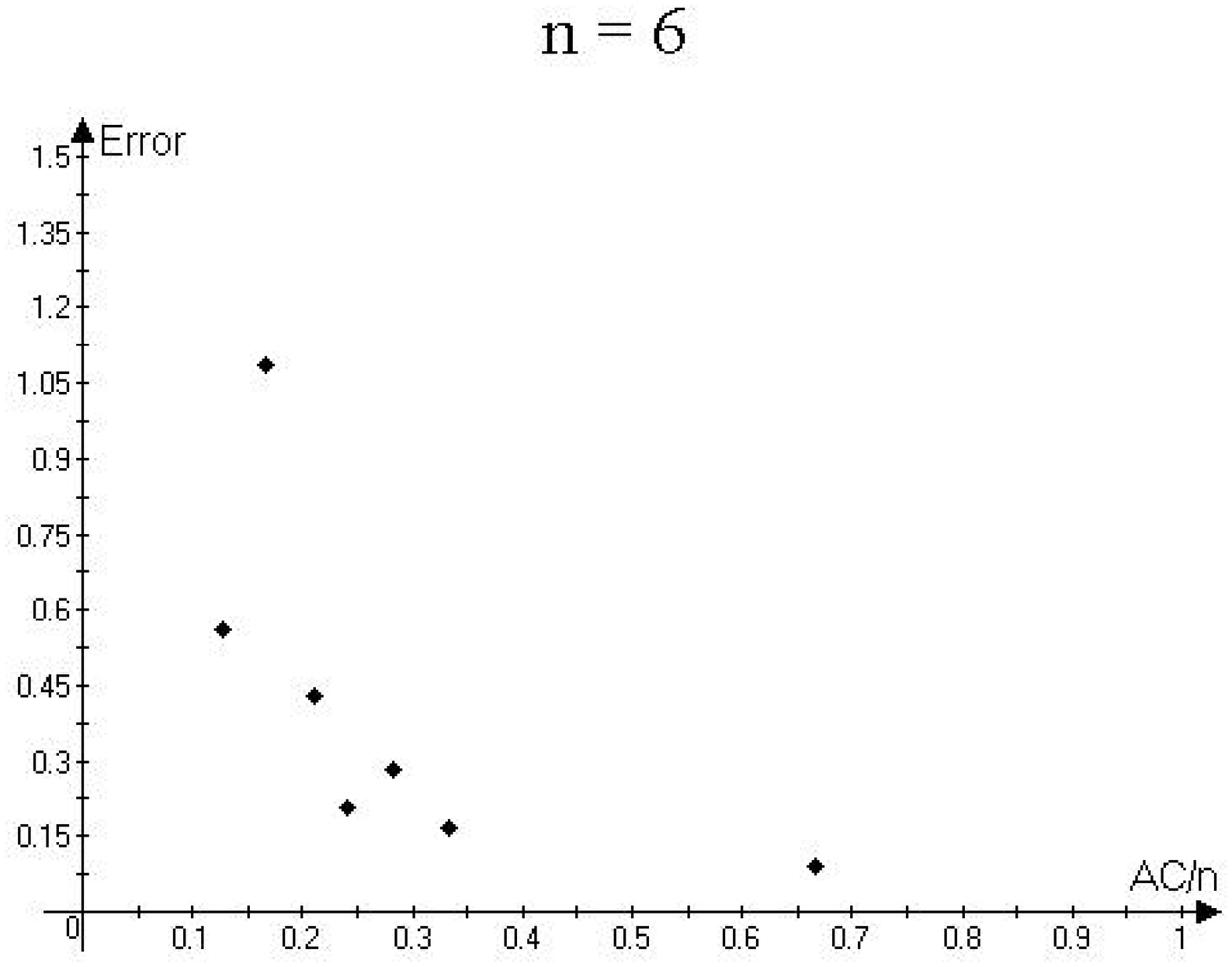} 
\hfill
\includegraphics[width=0.45\textwidth]{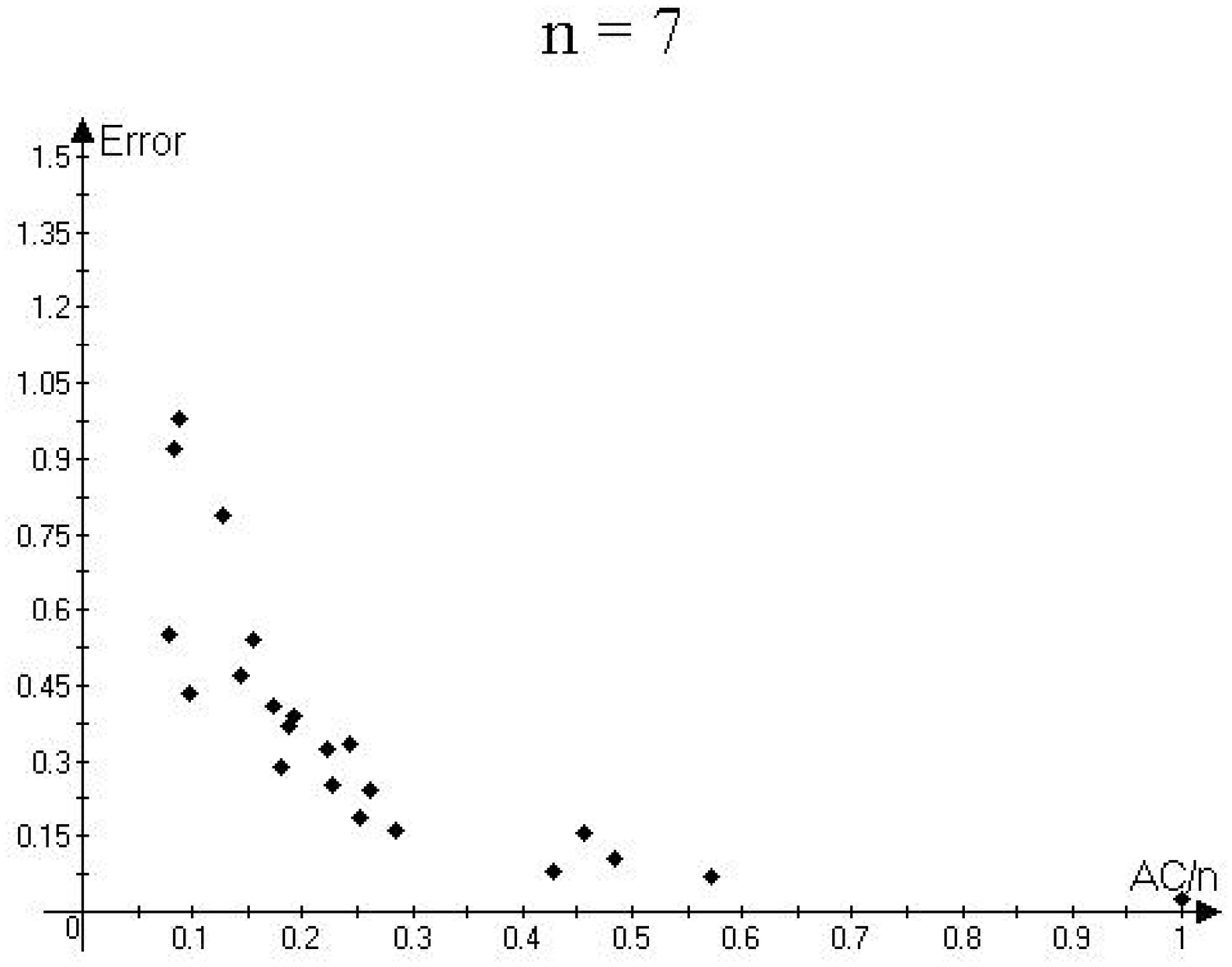} 
\end{figure}

\begin{figure}[h]
\vspace{0mm}
\includegraphics[width=0.45\textwidth]{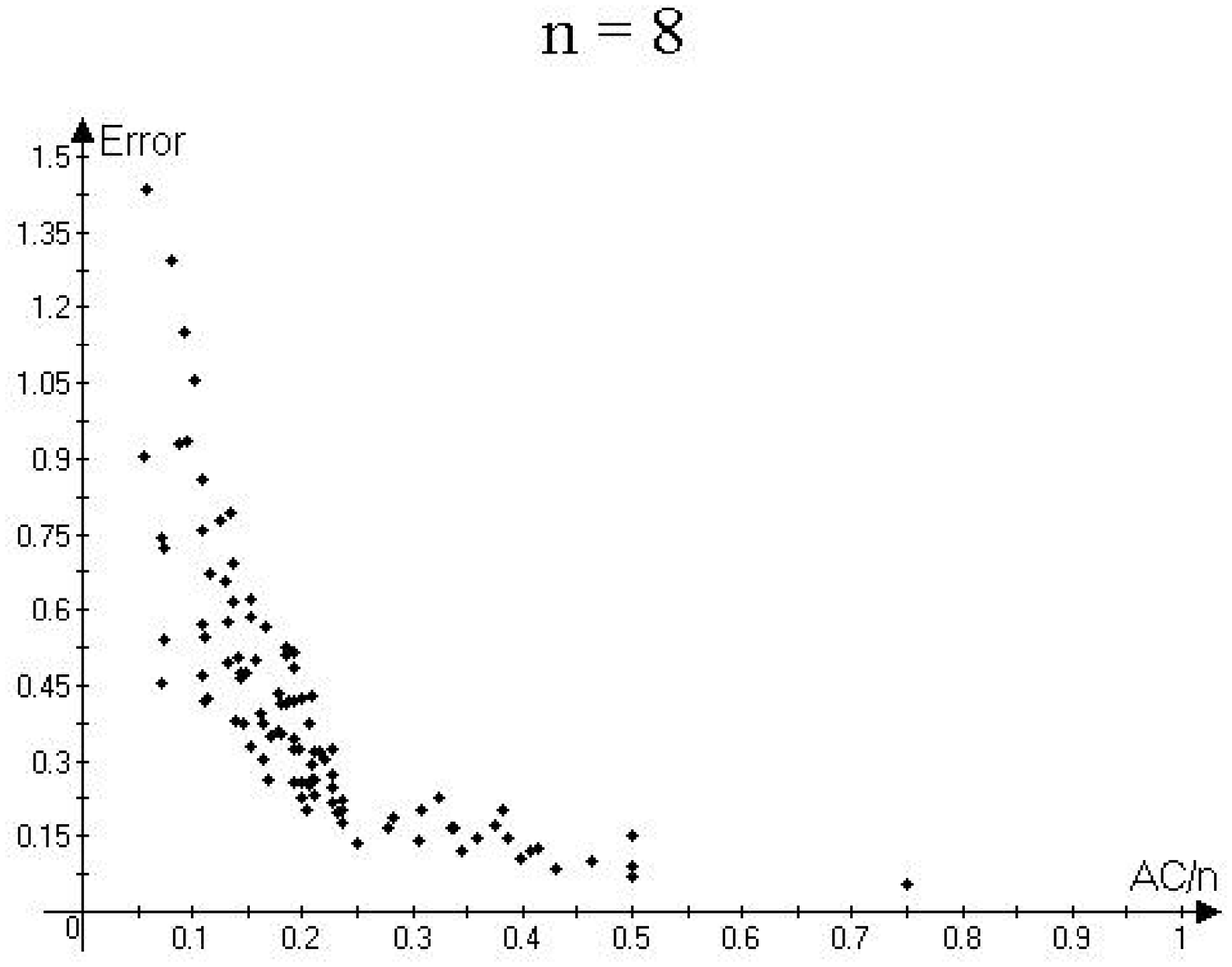} 
\hfill
\includegraphics[width=0.45\textwidth]{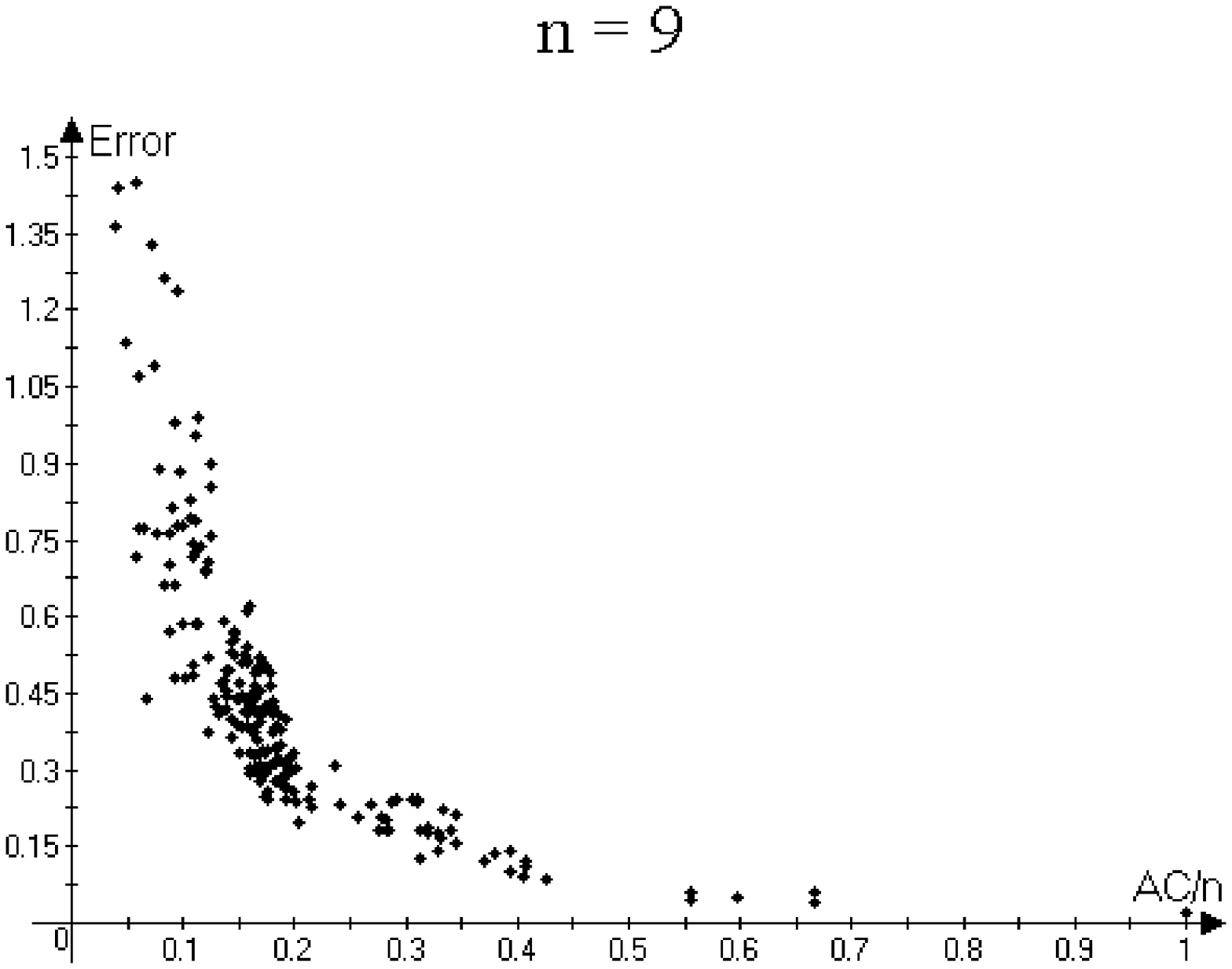} 
\end{figure}

\vspace{5mm}

The charts show, in particular, that  $Error$ decreases significantly with respect to the ratio $\lambda_2(G)/n$.

\subsection{Eulerian circuits}

To our knowledge (in contrast to the coundting of eulerian orientations) approximate polynomial algorithms for counting the number of Eulerian circuits 
have been obtained only in the literature fore some special classes of graphs having low density, see \cite{CCM2012} and \cite{TV2001}. 
However, we have the formula for $EC(G)$ similar to \eqref{main_eq} for $\gamma$-mixing graphs. This formula is more complicated, so we need some additional notations. Let
\begin{equation}
	W = \hat{Q}^{-1} = (Q + J)^{-1},
\end{equation}
where $Q$ is the Laplacian matrix and $J$ denotes the matrix with every entry $1$.
Let $\vec{\alpha} = (\alpha_1,\ldots,\alpha_n)\in \mathbb{R}^n$ be defined by
\begin{equation}
 	\alpha_j = W_{jj}
\end{equation}
Let $\vec{\beta} = Q\vec{\alpha}$ and
\begin{equation}\label{K_1}
	C_1 = \exp\left( - \sum_{j=1}^{n-1} \sum\limits_{k=j+1}\limits^{n} \beta_j W_{jk} \beta_k \right),
\end{equation} 
\begin{equation}\label{K_2}
	C_2 = \exp\left( -  \sum_{j=1}^{n} \frac{\beta_j^2}{2(d_j+1)} \right),
\end{equation}
where $d_j$ is the degree of $v_j\in VG$. Let 
\begin{equation}
 R(\vec{\theta}) = \mbox{tr} (\Lambda(\vec{\theta}) W \Lambda(\vec{\theta}) W),
\end{equation}
where $\mbox{tr}(\cdot)$ is the trace fucntion, $\Lambda(\vec{\theta})$ denotes the diagonal
matrix whose diagonal elements are equal to components of the vector $Q\vec{\theta}$.
Let $\vec{e}^{(k)} = (e^{(k)}_1,\ldots,e^{(k)}_n)  \in \mathbb{R}^n$ be defined by $e^{(k)}_j = \delta_{jk}$, where
$\delta_{jk}$ is the Kronecker delta. Let $r_k = R(\vec{e}^{(k)})$,
\begin{equation}\label{K_3}
	C_3 = \exp\left( \sum_{j=1}^{n} \frac{r_j}{2(d_j+1)} \right).
\end{equation}
Finally, let 
\begin{equation}\label{K_4}
	C_4 = \prod\limits_{\{v_j,v_k\}\in EG} P_{jk}, 
\end{equation}
where $P_{jk}$ is the same as in \eqref{main_eq}.

\begin{Proposition}\label{main1}
	Let $G$ be an undirected simple graph  with  $n$ vertices $v_1, v_2, \ldots, v_n$ having even degrees. Let $G$ be a $\gamma$-mixing graph 
	for some $\gamma>0$. 
	Then 
	\begin{equation}\label{main_eq1}
		\begin{aligned}
		EC(G) = \left(1 + \delta'(G) \right)
		\left(
		   \prod\limits_{j=1}\limits^{n}\left(\frac{d_j}{2}-1\right)! 
		   \,2^{|EG|-\frac{n-1}{2}} \pi^{-\frac{n-1}{2}} { \sqrt{t(G)}} C_1 C_2 C_3 C_4
		\right),
			\end{aligned}
	\end{equation}
	where  $C_1,C_2,C_3,C_4$ are defined according to \eqref{K_1}, \eqref{K_2}, \eqref{K_3}, \eqref{K_4},
respectively, $d_j$ is  the degree of vertex $v_j$,	  $t(G)$ is the number of spanning trees of  $G$ 
and for any  $\varepsilon>0$ 
	\begin{equation}
		|\delta'(G)| \leq C' n^{-1/2+\varepsilon},
	\end{equation}
	where constant $C' > 0$ depends only on $\gamma$ and $\varepsilon$.   
\end{Proposition}

\begin{Remark}
	One can obtain for the case of $G=K_n$ that 
\begin{equation}
	\begin{aligned}
		 C_1C_2C_3C_4 = 1 + O(n^{-1})
	\end{aligned}
\end{equation}
	Using Remark 5.2 and Stirling's formula for factorials, the result of Proposition \ref{main1} for this case can be 
	reduced to the result of \cite{Brendan1995} 
	on counting Eulerian circuits in the complete graph.
\end{Remark}

We plan to give the proof of Proposition \ref{main1} in a subsequent paper. 
In the present work we just compare the answers given by  formula \eqref{main_eq1} against the exact
values for small graphs. 
Let
\begin{equation}
	Error'(G) = \frac{EC_{approx}(G) - EC(G)}{EC(G)},
\end{equation} 
where $EC_{approx}(G)$ is taken according to right-hand side of \eqref{main_eq1}.
The following charts show the dependence of $Error'(G)$ on the ratio $\lambda_2(G)/n$, where $\lambda_2(G)$ 
is the algebraic connectivity and $n=|VG|=6,7,8,9$:

\vspace{5mm}

\begin{figure}[h] 
\vspace{0mm}
\includegraphics[width=0.45\textwidth]{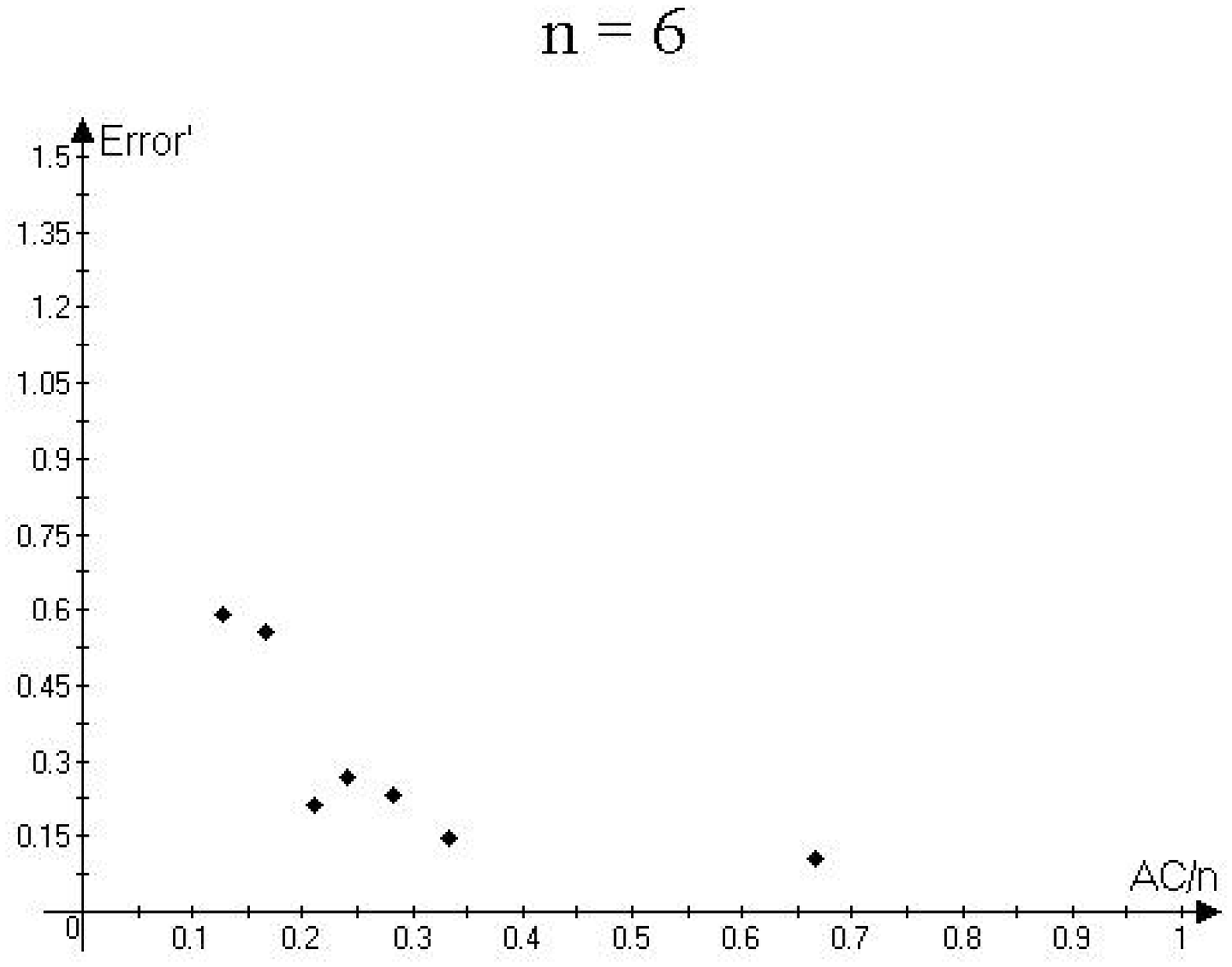}
\hfill
\includegraphics[width=0.45\textwidth]{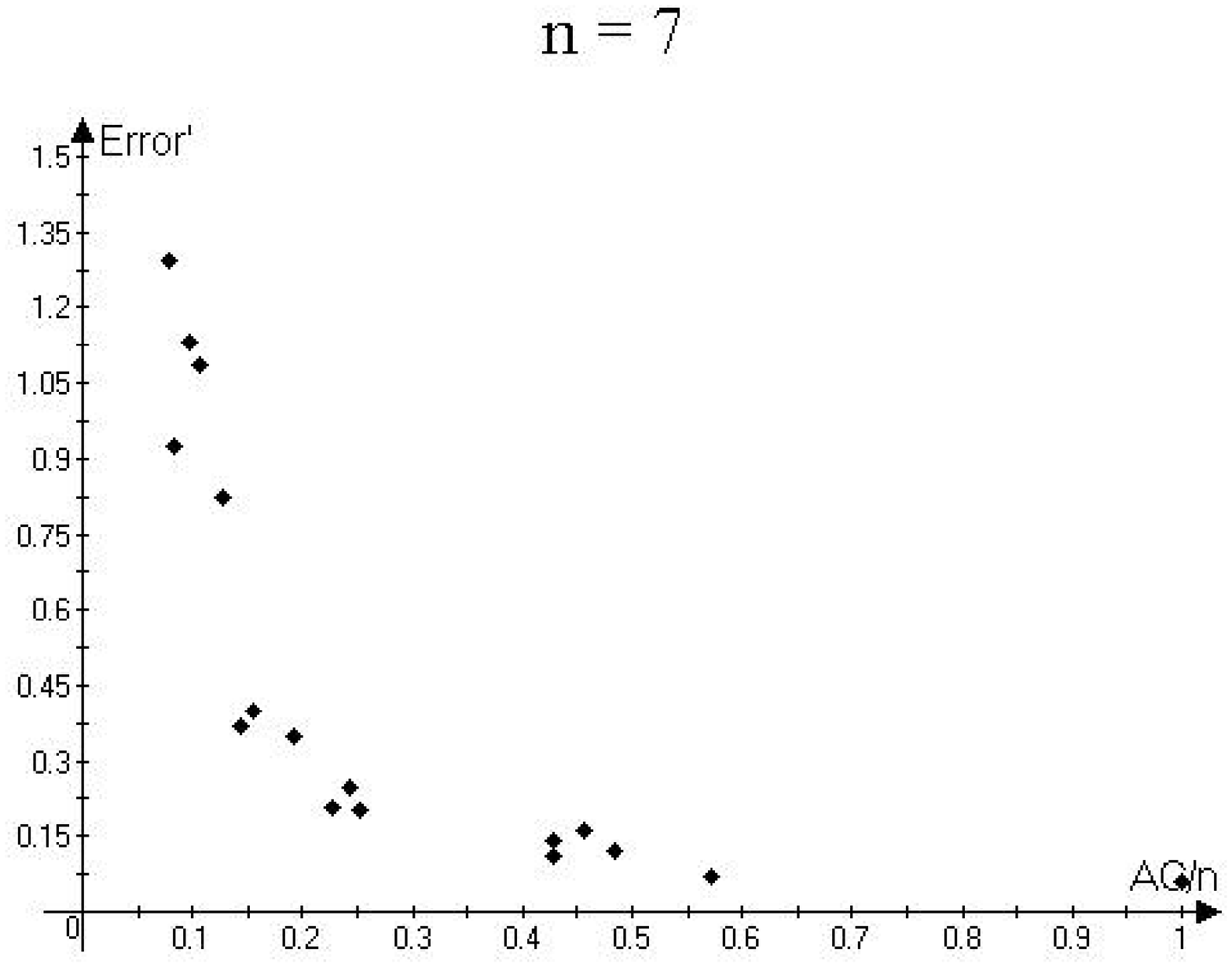} 
\end{figure}

\begin{figure}[h]

\includegraphics[width=0.45\textwidth]{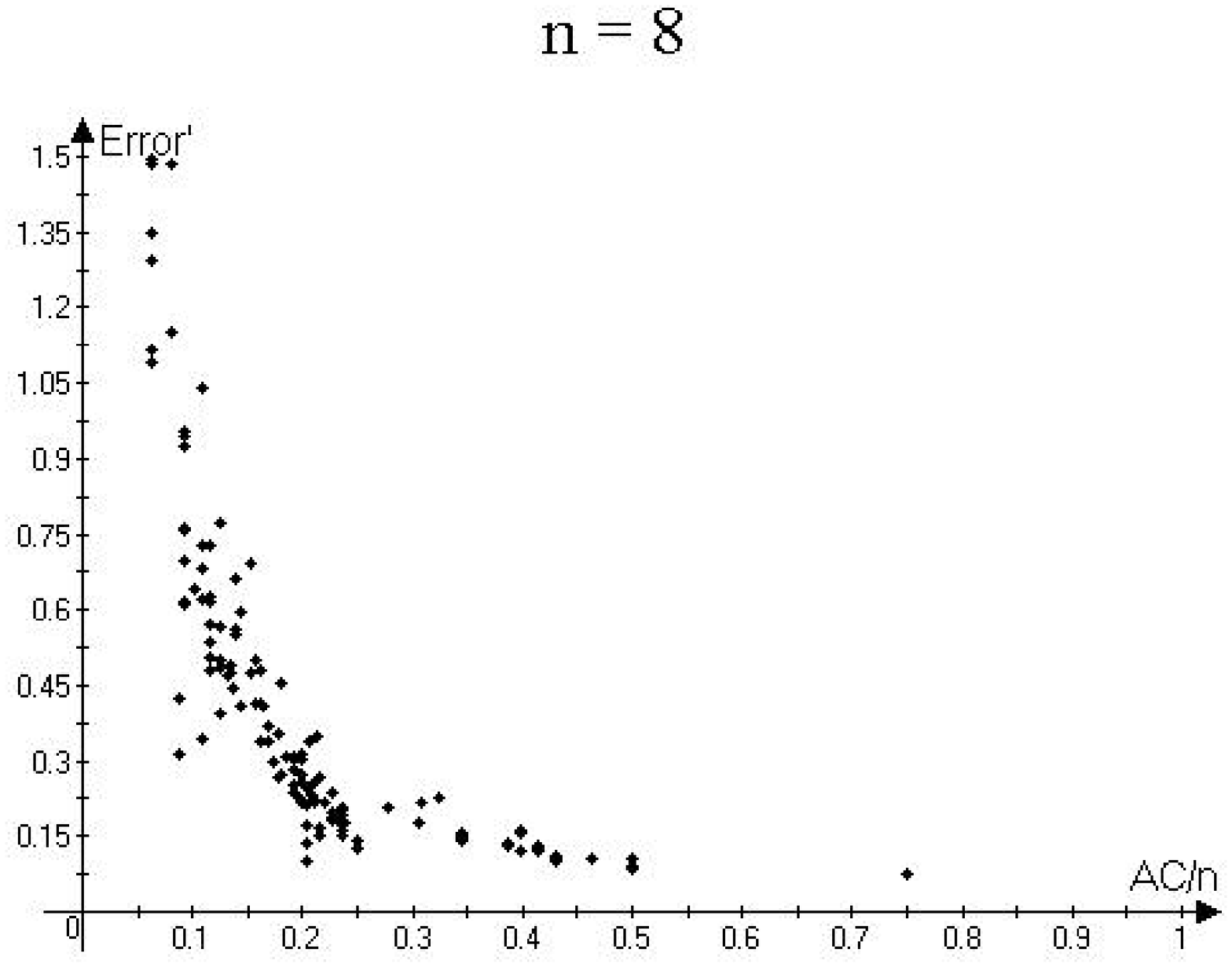}
\hfill
\includegraphics[width=0.45\textwidth]{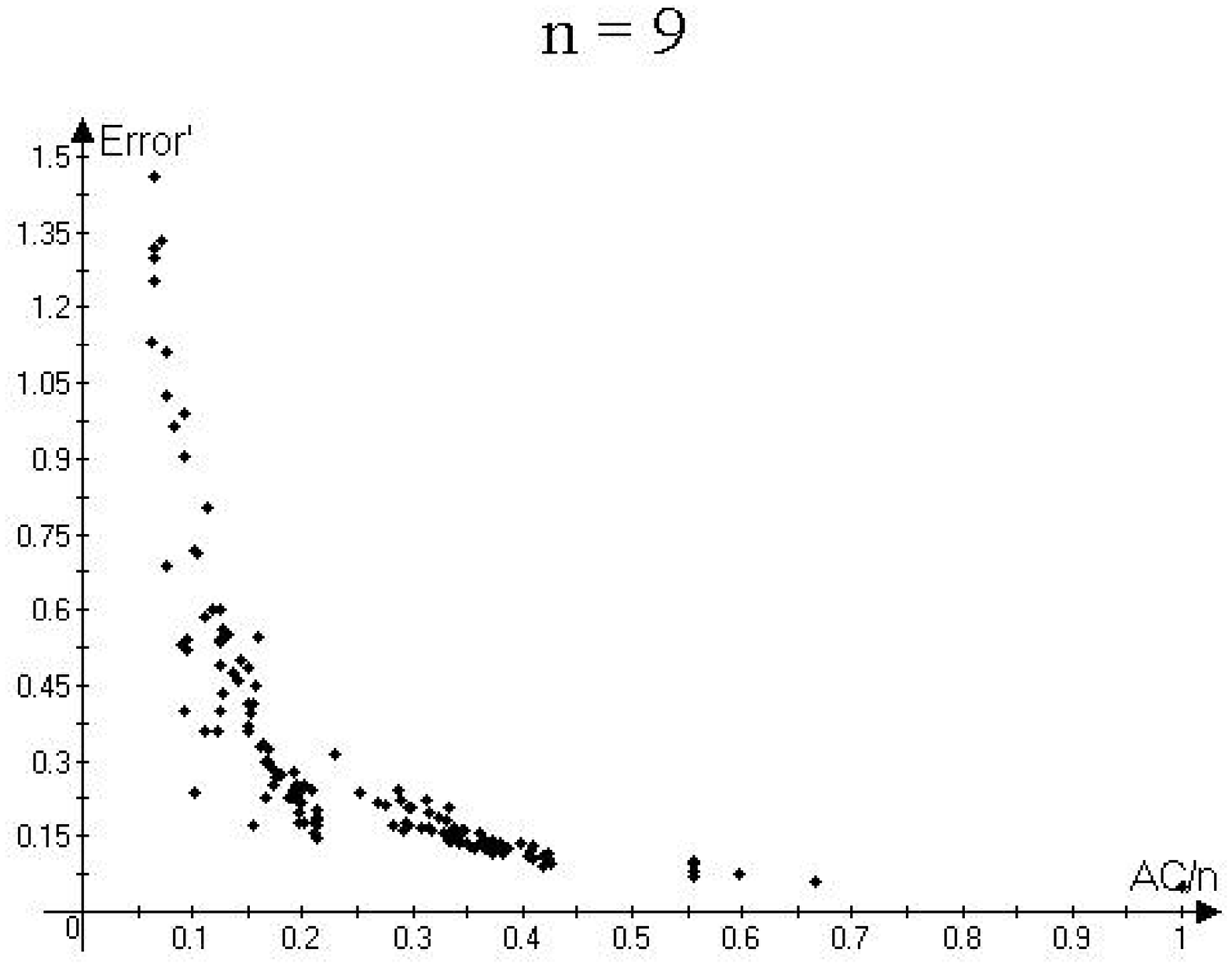} 
\end{figure}

\vspace{5mm}

The charts show, in particular, that  $Error'$ decreases significantly with respect to the ratio $\lambda_2(G)/n$.
\section*{Acknowledgements}
 This work was carried out under the supervision of S.P. Tarasov and supported in part by
RFBR grant no 11-01-00398a.

\noindent
{ {\bf M.I. Isaev}\\
Moscow Institute of Physics and Technology,

141700 Dolgoprudny, Russia\\
Centre de Math\'ematiques Appliqu\'ees, Ecole Polytechnique,

91128 Palaiseau, France\\
e-mail: \tt{Isaev.M.I@gmail.com}}\\

\noindent
{ {\bf K.V. Isaeva}\\
Moscow Institute of Physics and Technology,

141700 Dolgoprudny, Russia\\
e-mail: \tt{Isaeva.K.V@gmail.com}}

\end{document}